\documentclass[12pt]{article}

\title{\textbf{How real are real numbers?}}

\author{\textbf{Gregory Chaitin}\thanks{IBM T. J. Watson Research Center, P. O. Box 218,
Yorktown Heights, NY 10598, U.S.A., \emph{chaitin@us.ibm.com.}}}

\date{}

\begin{document}

\maketitle

\begin{abstract}
We discuss mathematical and physical arguments 
against continuity and in favor of discreteness,
with particular emphasis on the ideas of Emile Borel (1871--1956).
\end{abstract}

\section{Introduction}

Experimental physicists know how difficult accurate measurements are. 
No physical quantity has ever been measured with more than 15 or so digits of accuracy.
Mathematicians, however, freely fantasize with infinite-precision real numbers.
Nevertheless within pure math the notion of a real number is extremely problematic.

We'll compare and contrast 
two parallel historical episodes:
\begin{enumerate}
\item
the diagonal and probabilistic proofs that reals are uncountable, and
\item
the diagonal and probabilistic proofs that there are uncomputable reals.
\end{enumerate}
Both case histories open chasms beneath the feet of mathematicians.
In the first case these are the famous Jules Richard paradox (1905),  Emile Borel's know-it-all real
(1927), and the fact that most reals are unnameable, which was the subject of 
[Borel, 1952], his last book, published when Borel was 81 years old [James, 2002].

In the second case the frightening features are
the unsolvability of the halting problem (Turing, 1936), the fact that most reals are uncomputable,
and last but not least, the 
halting probability $\Omega$, which is irreducibly complex (algorithmically random),
maximally unknowable,
and dramatically illustrates the limits of reason
[Chaitin, 2005].

In addition to this mathematical soul-searching
regarding real numbers,
some physicists are beginning to 
suspect that the physical universe is actually discrete 
[Smolin, 2000]
and perhaps even a giant computer [Fredkin, 2004, Wolfram, 2002].
It will be interesting to see how far this 
so-called
``digital philosophy,'' ``digital physics'' viewpoint can be taken.

\emph{Nota bene:} To simplify matters, throughout this paper 
we restrict ourselves to reals in the interval between 0 and 1. 
We can therefore identify a real number
with the infinite sequence of digits or bits after its decimal or binary point.

\section{Reactions to Cantor's Theory of Sets: The Trauma of the Paradoxes of Set Theory}

Cantor's theory of infinite sets, developed in the late 1800's, was a decisive advance
for mathematics, but it provoked raging controversies and abounded in paradox.
One of the first books by the distinguished French mathematician Emile Borel (1871--1956)\footnote
{For a biography of Borel, see [James, 2002].} 
was his \emph{Le\c{c}ons sur la Th\'eorie des Fonctions} [Borel, 1950], 
originally published in 1898, and
subtitled \emph{Principes de la th\'eorie des ensembles en vue des applications \`a la th\'eorie
des fonctions.}

This was one of the first books promoting Cantor's theory of sets (\emph{ensembles}), 
but Borel had serious
reservations about certain aspects of Cantor's theory, 
which Borel kept adding to later editions of his book
as new appendices. The final version of Borel's book, 
which was published by Gauthier-Villars in 1950, has been kept in print
by Gabay. That's the one that I have, and
this book is a treasure trove of interesting mathematical, philosophical
and historical material.

One of Cantor's crucial ideas is the distinction between
the denumerable or countable infinite sets,
such as the positive integers or the rational numbers, 
and the much larger nondenumerable or uncountable infinite sets, 
such as the real numbers or the points in the plane or in space.
Borel had constructivist leanings, and as we shall see he felt comfortable with denumerable sets,
but very uncomfortable with nondenumerable ones.
And one of Cantor's key results that is discussed by Borel is Cantor's proof that the set of
reals is nondenumerable, i.e., cannot be placed in a one-to-one correspondence with the positive
integers.
I'll prove this now in two different ways.

\subsection{Cantor's diagonal argument: Reals are uncountable/nondenumerable}

Cantor's proof of this is a \emph{reductio ad absurdum.}

Suppose on the contrary that we have managed to list all the reals, with a first real, a second
real, etc.   Let $d(i,j)$ be the $j$th digit after the decimal point of the $i$th real in the list.
Consider the real $r$ between 0 and 1 whose $k$th digit 
is defined to be 4 if $d(k,k) = 3$, and 3 otherwise.
In other words, we form $r$ by taking all the decimal digits 
on the \textbf{diagonal} of the list of all reals,
and then changing each of these diagonal digits.

The real $r$ differs from the $i$th real in this presumably complete list of all reals, 
because their $i$th digits are different.
Therefore this list cannot be complete, and the set of reals is uncountable. \emph{Q.E.D.}

\emph{Nota bene:}
The most delicate point in this proof is to avoid having $r$ end in an infinity of 0's 
or an infinity of 9's,
to make sure that having its $k$th digit differ from the $k$th digit of the $k$th real
in the list suffices to guarantee that $r$ is not equal to the $k$th real in the list.
This is how we get around the fact that some reals can have more than one decimal representation.

\subsection{Alternate proof: Any countable/denumerable set of reals has measure zero}

Now here is a radically different proof that the reals are uncountable.
This proof, which I learned in [Courant \& Robbins, 1947], was perhaps or at least 
\textbf{could have been} originally 
discovered by Borel, because it uses the mathematical notion of \emph{measure},
which was invented by Borel and later perfected by his 
Ecole Normale Sup\'erieure student Lebesgue,
who now usually gets all the credit.

Measure theory and probability theory are really one and the same---it's just 
different names for the same concepts.
And Borel was interested
in both the technical mathematical aspects 
and in the many important practical
applications, which Borel discussed in many of his books.

So let's suppose we are given a real $\epsilon>0$, which we shall later make arbitrarily small.
Consider again that supposedly complete enumeration of all the reals, a first one, a second one, etc.
Cover each real with an interval, and take the interval for covering the $i$th real in the list 
to be of length $\epsilon/2^i$.  The total length of all the covering intervals is therefore
\[
   \frac{\epsilon}{2} + \frac{\epsilon}{4} + \cdots \frac{\epsilon}{2^i} + \cdots 
   \;\;\; = \;\;\; \epsilon,
\]
which we can make as small as we wish.

In other words, any countable set of reals has \emph{measure zero} and is a so-called \emph{null set,}
i.e., has zero probability and is an infinitesimal subset of the set of all reals.
\emph{Q.E.D.}

We have now seen
the two fundamentally different ways of showing that the reals are infinitely more
numerous than the positive integers, i.e., that the set of all reals is a higher-order
infinity than the set of all positive integers.

So far, so good!
But now, let's show what a minefield this is.

\subsection{Richard's paradox: Diagonalize over all nameable reals
$\longrightarrow$ a nameable, unnameable real}

The problem is that the set of reals is uncountable, 
but the set of all possible texts in English or French is countable,
and so is the set of all possible mathematical definitions or the set of all possible
mathematical questions, since these also have to be formulated within a language, yielding
at most a denumerable infinity of possibilities.
So there are too many reals, and not enough texts.

The first person to notice this difficulty was Jules Richard in 1905, and 
the manner in which he formulated the problem  
is now called Richard's paradox.

Here is how it goes. Since all possible texts in French (Richard was French) can be listed or
enumerated, a first text, a second one, etc.,\footnote
{List all possible texts in size order, 
and within texts that are the same size, in alphabetical order.}
you can diagonalize over all the reals that
can be defined or named in French and produce a real number that cannot be defined 
and is therefore unnameable.
However, we've just indicated how to define it or name it!

In other words, Richard's paradoxical real differs from every real that is definable in French,
but nevertheless can itself be defined in French by specifying in detail how to apply Cantor's
diagonal method to the list of all possible mathematical definitions for individual real numbers
in French!

How very embarrassing!  Here is a real number that is simultaneously nameable yet at the same time
it cannot be named using any text in French.

\subsection{Borel's know-it-all number}

The idea of being able to list or enumerate all possible texts in a language is an extremely
powerful one, and it was exploited by Borel in 1927 [Tasi\'c, 2001, Borel, 1950]
in order to define a real number that can answer
every possible yes/no question!

You simply write this real in binary, and use the $n$th bit of its binary expansion to answer
the $n$th question in French.

Borel speaks about this real number ironically. 
He insinuates that it's illegitimate, unnatural, artificial,
and that it's an ``unreal'' real number, one that there is no reason
to believe in.

Richard's paradox and Borel's number are discussed in [Borel, 1950] on the pages given in the
list of references, but the next paradox was considered so important by Borel that he devoted
an entire book to it.  In fact, this was Borel's last book [Borel, 1952] and it was published,
as I said,
when Borel was 81 years old.  I think that when Borel wrote this work
he must have been thinking about
his legacy, since this was to be his final book-length mathematical statement.  
The Chinese, I believe, place special value on
an artist's final work, considering that in some sense it contains or
captures that artist's soul.\footnote
{I certainly feel that way about Bach's \emph{Die Kunst der Fuge} and about Bergman's
\emph{Fanny och Alexander}.}
If so, [Borel, 1952] is Borel's ``soul work.''

Unfortunately I have not been able to obtain this crucial book. 
But based on a number of remarks
by other people and based on what I do know about Borel's methods and concerns, I am
fairly confident that I know what [Borel, 1952] contains.
Here it is:

\subsection{Borel's ``inaccessible numbers:'' Most reals are unnameable, with probability one}

Borel's often-expressed credo is that a real number is really real 
only if it can be expressed, only if it can be uniquely defined, using a finite
number of words.\footnote{See for example [Borel, 1960].} 
It's only real 
if it can be named or specified as an individual mathematical object.
And in order to do this we must necessarily employ some particular language, e.g., French.
Whatever the choice of language, there will only be a countable infinity of possible texts,
since these can be listed in size order, and among texts of the same size, in alphabetical order.

This has the devastating consequence that there are only a denumerable
infinitely of such ``accessible'' reals, and therefore, as we saw in Sec.\ 2.2,
the set of accessible reals has measure zero.

So, in Borel's view, most reals, with probability one, are mathematical fantasies, because
there is no way to specify them uniquely.  Most reals are inaccessible to us,
and will never, ever, be picked out as individuals using \textbf{any} conceivable mathematical tool,
because whatever these tools may be they could always be explained in French, and therefore can only
``individualize''
a countable infinity of reals, a set of reals of measure zero,
an infinitesimal subset of the set of all possible reals. 

Pick a real at random, and the probability is zero that it's accessible---the probability is zero
that it \textbf{will ever be} accessible to us as an individual mathematical object.

\section{History Repeats Itself: Computability Theory and Its Limitative Meta-Theorems}

That was an exciting chapter in the history of ideas, wasn't it!
But history moves on, and the collective attention of the human species shifts elsewhere,
like a person who is examining a huge painting.

What completely transformed the situation is the
\textbf{idea} of the computer, the computer as a mathematical concept, not a practical device,
although the current ubiquity of computers doesn't hurt.  It is, as usual, unfair to
single out an individual, but in my opinion the crucial event was the 1936 paper by
Turing \emph{On computable numbers,} and here Turing is in fact referring to computable real numbers.
You can find this paper at the beginning of the collection [Copeland, 2004],  
and at the end of this book there happens to be 
a much more understandable paper by Turing explaining just the key idea.\footnote
{It's Turing's 1954 Penguin \emph{Science News} paper on \emph{Solvable and unsolvable problems,}
which I copied out into a notebook by hand when I was a teenager.}

History now repeats itself and recycles the ideas that were presented in Sec.\ 2. 
This time the texts will be written in artificial formal languages, they 
will be computer programs or proofs in a formal axiomatic math theory. 
They won't be texts that are written in a natural language like English or French.
And this time we won't get paradoxes, instead we'll get meta-theorems, we'll get
limitative theorems, ones that show the limits of computation or 
the limitations of formal math theories.
So in their current reincarnation, which we'll now present,
the ideas that we saw in Sec.\ 2 definitely become much sharper and clearer.

Formal languages avoid the paradoxes by removing
the ambiguities of natural languages. The paradoxes are eliminated, but there is a price.
Paradoxical natural languages are evolving open systems. Artificial languages are
static closed systems subject to limitative meta-theorems.
You avoid the paradoxes, but you are left with a corpse!

The following tableau summarizes the transformation (paradigm shift):
{\footnotesize
\begin{itemize}
\item
Natural languages $\longrightarrow$ Formal languages.
\item
Something is true $\longrightarrow$ 
Something is provable within a particular formal axiomatic math theory.\footnote
{This part of the paradigm shift is particularly important in
the story of how G\"odel converted the paradox of ``this statement is false''
into the proof of his famous 1931 incompleteness theorem, which is based on
``this statement is unprovable.''
This changes something that's true if and only if it's false, into something that's
true if and only if it's unprovable, thus transforming a paradox into a meta-theorem.}
\item
Naming a real number $\longrightarrow$ Computing a real number digit by digit.
\item
Number of words required to name something\footnote{See [Borel, 1960].} $\longrightarrow$ 
Size in bits of the smallest program for computing something 
(program-size complexity).\footnote{See [Chaitin, 2005].}
\item
List of all possible texts in French $\longrightarrow$ List of all possible programs, or 
\\
List of all possible texts in French $\longrightarrow$ List of all possible proofs.\footnote
{The idea of systematically combining concepts in every possible way can be traced
through Leibniz back to Ramon Llull (13th century), and is ridiculed by Swift in
\emph{Gulliver's Travels} (Part III, Chapter 5, on the Academy of Lagado).}
\item
Paradoxes $\longrightarrow$ Limitative meta-theorems.
\end{itemize}
}

Now let's do Sec.\ 2 all over again.
First we'll examine two different proofs that there are uncomputable
reals: a diagonal argument proof, and a measure-theoretic proof.
Then we'll show how the Richard paradox yields the unsolvability of the halting problem.
Finally we'll discuss 
the halting probability $\Omega$, which plays roughly the same role here that
Borel's know-it-all real did in Sec.\ 2.

\subsection{Turing diagonalizes over all computable reals
$\longrightarrow$ uncomputable real}

The set of all possible computer programs is countable, therefore the set of all
computable reals is countable, and diagonalizing over the computable reals immediately
yields an uncomputable real. \emph{Q.E.D.}

Let's do it again more carefully.

Make a list of all possible computer programs.
Order the programs by their size,
and within those of the same size, order them alphabetically.
The easiest thing to do is to include all the possible character strings
that can be formed from
the finite alphabet of the programming language, even though
most of these will be syntactically invalid programs.

Here's how we define the uncomputable diagonal number $0<r<1$.
Consider the $k$th program in our list. If it is syntactically invalid, or if
the $k$th program never outputs a $k$th digit,
or if the $k$th digit output by the $k$th program isn't a 3,
pick 3 as the $k$th digit of $r$.
Otherwise, if the $k$th digit output by the $k$th program is a 3, pick 4 as the $k$th digit of $r$.

This $r$ cannot be computable, because its $k$th digit is different from the $k$th digit
of the real number that is computed by the $k$th program, if there is one.
Therefore there are uncomputable reals, real numbers that cannot be calculated digit by digit
by any computer program.

\subsection{Alternate proof: Reals are uncomputable with probability one}

In a nutshell, the set of computer programs is countable, therefore the set of
all computable reals is countable, and therefore, as in Sec.\ 2.2, of measure zero.
\emph{Q.E.D.}

More slowly, consider the $k$th computer program again.  If it is syntactically invalid
or fails to compute a real number, let's skip it.  If it does compute a real,
cover that real with an interval of length $\epsilon/2^k$.  Then the total length
of the covering is less than $\epsilon$, which can be made arbitrarily small,
and the computable reals are a null set.

In other words, the probability of a real's being computable is zero, and the probability
that it's uncomputable is one.\footnote
{Who should be credited for this measure-theoretic proof that there are uncomputable reals?
I have no idea.  It seems to have always been part of my mental baggage.}

What if we allow arbitrary, highly nonconstructive means 
to specify particular reals, not just computer programs?
The argument of Sec.\ 2.5 carries over immediately  within our new framework
in which we consider formal languages instead of natural languages.
Most reals remain unnameable, with probability one.\footnote
{This theorem is featured in [Chaitin, 2005] at the end of the chapter entitled
\emph{The Labyrinth of the Continuum.}}

\subsection{Turing's halting problem: No algorithm settles halting, no 
formal axiomatic math theory settles halting}

Richard's paradox names an unnameable real.
More precisely, it diagonalizes over all reals uniquely specified by French texts
to produce a French text specifying an unspecifiable real.
What becomes of this in our new context in which we name reals by computing them?

Let's go back to Turing's use of the diagonal argument in Sec.\ 3.1.
In Sec.\ 3.1 we constructed an uncomputable real $r$.
It must be uncomputable, by construction.
Nevertheless, as was the case in the Richard paradox, it would seem that we gave
a procedure for calculating Turing's diagonal real $r$ digit by digit.
How can this procedure fail?
What could possibly go wrong?

The answer is this:  The only noncomputable step has got to be determining if the
$k$th computer program will \textbf{ever} output a $k$th digit.
If we could do that, then we could certainly compute the uncomputable real $r$ of Sec.\ 3.1.

In other words, Sec.\ 3.1 actually \textbf{proves} that there can be no algorithm for
deciding if the $k$th computer program will ever output a $k$th digit.

And this is a special case of what's called Turing's halting problem.
In this particular case, 
the question is whether or not the wait for a $k$th digit will ever terminate.
In the general case, the question is whether or not a computer program will ever halt.

The algorithmic unsolvability of Turing's halting problem is an extremely fundamental
meta-theorem.  
It's a much stronger result than G\"odel's famous 1931 incompleteness theorem.
Why? Because
in Turing's original 1936 paper he immediately points out 
how to derive incompleteness from the halting problem.

A formal axiomatic math theory (FAMT) consists of a finite set of axioms and of
a finite set of rules of inference for deducing the consequences of those axioms. 
Viewed from a great distance,
all that counts is that there is an algorithm for enumerating (or generating) all
the possible theorems, all the possible consequences of the axioms, one by one,
by systematically applying the rules of inference in every possible way.
This is in fact what's called a breadth-first (rather than a depth-first) tree walk,
the tree being the tree of all possible deductions.\footnote
{This is another way to achieve the effect of running through all possible texts.}

So, argued Turing in 1936, if there were a FAMT that always enabled you to decide whether or not
a program eventually halts, there would in fact be an algorithm for doing so.
You'd just run through all possible proofs until you find a proof that the program halts
or you find a proof that it never halts.

So uncomputability is much more fundamental than incompleteness.
Incompleteness is an immediate corollary of uncomputability.
But uncomputability is \textbf{not} a corollary of incompleteness.
The concept of incompleteness does not contain the concept of uncomputability.

Now let's get an even more disturbing limitative meta-theorem.
We'll do that by considering the halting probability $\Omega$ [Chaitin, 2005], which 
is what corresponds to
Borel's know-it-all real (Sec.\ 2.4) in the current context.\footnote
{[Tasi\'c, 2001] was the first person to make the connection between Borel's real and $\Omega$.
I became aware of Borel's real through Tasi\'c.}

\subsection{Irreducible complexity, perfect randomness, maximal unknowability: 
The halting probability $\Omega$}

Where does the halting probability come from?  Well, our motivation is the contrast between
Sec.\ 3.1 and Sec.\ 3.2. Sec.\ 3.1 is to Sec.\ 3.2 as the halting problem is to the halting
probability!  In other words, the fact that we found an easier way to show the existence of
uncomputable reals using a probabilistic argument,
suggests looking at the probability that a program chosen at random will ever halt instead
of considering individual programs as in Turing's 1936 paper.

Formally, the halting probability $\Omega$ is defined as follows:
\[
0 \;\;\; < \;\;\; \Omega \;\;\; \equiv \;\;\;
\sum_{\mbox{\scriptsize program $p$ halts}} 2^{-(\mbox{\scriptsize the size in bits of $p$})}
\;\;\; < \;\;\; 1.
\]  
To avoid having this sum diverge to infinity instead of converging to a number between zero and one,
it is important that the programs $p$ should be self-delimiting (no extension of a valid program
is a valid program; see [Chaitin, 2005]).

What's interesting about $\Omega$ is that it behaves like a compressed version of Borel's
know-it-all real.  Knowing the first $n$ bits of Borel's real enables us to answer the
first $n$ yes/no questions in French.
Knowing the first $n$ bits of $\Omega$ enables us to answer the halting problem for all
programs $p$ up to $n$ bits in size. 
I.e., $n$ bits of $\Omega$ tells us whether or not each $p$ up to $n$ bits in size ever halts.
(Can you see how?)
That's a lot of information!

In fact, $\Omega$ compactly encodes 
so much information that you essentially need an $n$-bit FAMT in order to be able to determine 
$n$ bits of $\Omega$!  In other words, $\Omega$ is \textbf{irreducible mathematical information},
it's a place where reasoning is completely impotent. The bits of $\Omega$ are mathematical 
facts that can be proved, but essentially only by adding them one by one as new axioms!
I'm talking about how difficult it is to prove theorems such as
\[
   \mbox{``the 5th bit of $\Omega$ is a 0''}
\]
and
\[
   \mbox{``the 9th bit of $\Omega$ is a 1''}
\]
or whatever the case may be.

To prove that $\Omega$ is computationally and therefore logically irreducible,
requires a theory of program-size complexity that I call
algorithmic information theory (AIT) [Chaitin, 2005].
The key idea in AIT is to measure the complexity of something via the size in bits
of the smallest program for calculating it.  This is a more refined version of Borel's idea
[Borel, 1960] of defining the complexity of a real number to be the number of words required
to name it. 

And the key fact that is proved in AIT about $\Omega$ is that 
\[
   H(\Omega_n) \geq n-c.
\]
I.e., 
\[
   \mbox{(the string $\Omega_n$ consisting of the first $n$ bits of $\Omega$)}
\]
has program-size complexity or ``algorithmic entropy $H$'' greater than or equal to $n-c$. 
Here $c$ is a constant,
and I'm talking about the size in bits of self-delimiting programs.

In other words, any self-delimiting program for computing the first $n$ bits of $\Omega$
will have to be at least $n-c$ bits long.

The irreducible sequence of bits of $\Omega$ is 
a place where mathematical truth has absolutely no pattern or structure
that we will ever be able to detect.
It's a place where mathematical truth has maximum
possible entropy---a place where, in a sense, God plays dice.\footnote
{On the other hand, if G\"odel is correct in thinking that mathematical intuition
can at times directly perceive the Platonic world of mathematical ideas,
then the bits of $\Omega$ may in fact be accessible.}

Why should we believe in real numbers, if most of them are uncomputable? 
Why should we believe in real numbers,  if most of them, it turns out,\footnote
{See the chapter entitled \emph{The Labyrinth of the Continuum} in [Chaitin, 2005].}
are maximally unknowable like $\Omega$?\footnote
{In spite of the fact that most individual real numbers will forever escape us, 
the notion of an arbitrary real has beautiful mathematical
properties and is a concept that helps us to organize and understand the real world.
Individual concepts in a theory \textbf{do not} need to have concrete meaning on their own; 
it is enough if the theory \textbf{as a whole} can be compared with the results of experiments.}

\section{Digital Philosophy and Digital Physics}

So much for mathematics! Now let's turn to physics.

Discreteness entered modern science through chemistry, when it was discovered that 
matter is built up out of atoms and molecules.  
Recall that the first experimental evidence for this was
Gay-Lussac's discovery of the simple integer ratios 
between the volumes of gaseous substances that are combined in
chemical reactions. This was the first evidence, two centuries ago,
that discreteness plays an important role in the physical world.

At first it might seem that quantum mechanics (QM), 
which began with Einstein's photon as the explanation for
the photoelectric effect in 1905, goes further in the direction of discreteness.
But the wave-particle duality discovered by de Broglie in 1925 is at the heart of QM,
which means that this theory is profoundly ambiguous regarding the question of discreteness vs.\
continuity.  QM can have its cake and eat it too, because discreteness 
is modeled via standing waves (eigenfunctions) in a continuous medium.

The latest strong hints in the direction of discreteness come from quantum gravity [Smolin, 2000],
in particular from the Bekenstein bound and the so-called ``holographic principle.''
According to these ideas the amount of information in any physical system is bounded,
i.e., is a finite number of 0/1 bits.

But it is not just fundamental physics that is pushing us in this direction.
Other hints come from our pervasive digital technology, from molecular biology
where DNA is the digital software for life, and from \emph{a priori}
philosophical prejudices going back to the ancient Greeks.

According to Pythagoras everything is number, and God is a mathematician.
This point of view has worked pretty well throughout the development of modern science.
However now a neo-Pythagorian doctrine is emerging, according to which everything is 0/1 bits, 
and the world is built entirely out of digital information.
In other words, now 
everything is software, God is a computer programmer, not a mathematician, 
and the world is a giant information-processing
system, a giant computer
[Fredkin, 2004, Wolfram, 2002, Chaitin, 2005].

Indeed,
the most important thing in understanding a complex system is to understand how it processes
information. This viewpoint regards physical systems as information processors,
as performing computations.  
This approach also sheds new light on microscopic quantum systems, as is 
demonstrated in the highly developed field of quantum information and quantum computation. 
An extreme version of this doctrine would attempt to build the world
entirely out of discrete digital information, out of 0 and 1 bits.\footnote
{This idea, like so many others, can be traced back to Leibniz.
He thought it was important enough to have it cast in the form of a medallion.}

Whether or not this ambitious new research program can eventually succeed, 
it will be interesting to see how far it gets.
The problem of the infinite divisibility of space and time has been with us for more than 
two millennia,
since Zeno of Elea and his famous paradoxes, 
and it is also discussed by Maimonides in his \emph{Guide for the Perplexed}
(12th century).

Modern versions of this ancient problem are, for example, the infinite amount of energy 
contained in the electric field surrounding a point electron
according to Maxwell's theory of electromagnetism, 
and the breakdown of space-time because of the formation of black holes due to extreme
quantum fluctuations (arbitrarily high energy virtual pairs) in the vacuum quantum field.

I do not expect that the tension between the continuous and the discrete 
will be resolved any time soon.
Nevertheless, one must try.
And, as we have seen in our two case studies, 
before being swept away,
each generation contributes something to the ongoing discussion.

\section*{References}

\footnotesize
Borel, E. [1950] \emph{Le\c{c}ons sur la Th\'eorie des Fonctions} (Gabay, Paris) pp.\ 161, 275.
\\
Borel, E. [1952] \emph{Les Nombres Inaccessibles} (Gauthier-Villars, Paris).
\\
Borel, E. [1960] \emph{Space and Time} (Dover, Mineola) pp.\ 212--214.
\\
Chaitin, G. [2005] \emph{Meta Math!}\ (Pantheon, New York),
math.HO/0404335, math.HO/0411091.
\\
Copeland, B. J. [2004] \emph{The Essential Turing} (Clarendon Press, Oxford).
\\
Courant, R. \& Robbins, H. [1947] \emph{What is Mathematics?}
(Oxford University Press, New York) Sec.\ II.4.2, pp.\ 79--83.
\\
Fredkin, E. [2004] \emph{http://www.digitalphilosophy.org}
\\
James, I. [2002] \emph{Remarkable Mathematicians} (Cambridge University Press, Cambridge)
pp.\ 283--292.
\\
Smolin, L. [2000] \emph{Three Roads to Quantum Gravity} (Weidenfeld \& Nicolson, London)
Chap.\ 8, pp.\ 95--105, Chap.\ 12, pp.\ 169--178.
\\
Tasi\'c, V. [2001] \emph{Mathematics and the Roots of Postmodern Thought} 
(Oxford University Press, New York) pp.\ 52, 81--82.
\\
Wolfram, S. [2002] \emph{A New Kind of Science} (Wolfram Media, Champaign).

\end{document}